\documentclass[11pt]{article}
\usepackage{amssymb}
\usepackage{verbatim,cite}
\usepackage{latexsym}
\usepackage{amsmath}
\usepackage{eucal}

\newtheorem{sub}{\name}[section]
\newcommand{\bs}{
\begin{sub}}
\newcommand{\es}{
\end{sub}}
\newcommand{\bsl}[1]{
\begin{sub}\label{#1}}
\newcommand{\bth}[1]{\def\name{Theorem}
\begin{sub}\label{t:#1}}
\newcommand{\blemma}[1]{\def\name{Lemma}
\begin{sub}\label{l:#1}}
\newcommand{\bcor}[1]{\def\name{Corollary}
\begin{sub}\label{c:#1}}
\newcommand{\bdef}[1]{\def\name{Definition}
\begin{sub}\label{d:#1}}
\newcommand{\bprop}[1]{\def\name{Proposition}
\begin{sub}\label{p:#1}}
\newcommand{\brem}[1]{\def\name{Remark}
\begin{sub}\label{r:#1}}
\newcommand{\bex}[1]{\def\name{Example}
\begin{sub}\label{e:#1}}

\newcommand{\rprop}[1]{Proposition~\ref{p:#1}}


\newcommand{\BA}{
\begin{array}}
\newcommand{\EA}{
\end{array}}

\newcommand{\BAN}{\renewcommand{\arraystretch}{1.2}
\setlength{\arraycolsep}{2pt}
\begin{array}}
\newcommand{\BAV}[2]{\renewcommand{\arraystretch}{#1}
\setlength{\arraycolsep}{#2}
\begin{array}}
\newcommand{\BSA}{
\begin{subarray}}
\newcommand{\ESA}{
\end{subarray}}
\newcommand{\BAL}{
\begin{aligned}}
\newcommand{\EAL}{
\end{aligned}}
\newcommand{\BALG}{
\begin{alignat}}
\newcommand{\EALG}{
\end{alignat}}
\newcommand{\BALGN}{
\begin{alignat*}}
\newcommand{\EALGN}{
\end{alignat*}}
\newcommand{\note}[1]{\textit{#1.}\hspace{2mm}}
\newcommand{\Proof}{\note{Proof}}
\newcommand{\qeda}{\hspace{10mm}\hfill $\square$}

\newcommand{\Remark}{\note{Remark}}

\newcommand{\abs}[1]{\left |#1\right |}
\newcommand{\norm}[1]{\left \|#1\right \|}

\def\angb<#1>{\langle #1 \rangle}

\newcommand{\opname}[1]{\mbox{\rm #1}\,}

\newcommand{\dist}{\opname{dist}}
\newcommand{\myfrac}[2]{{\displaystyle \frac{#1}{#2} }}
\newcommand{\myint}[2]{{\displaystyle \int_{#1}^{#2}}}
\newcommand {\dint}{{\displaystyle \int\!\!\int}}

\newcommand{\prt}{
\partial}

\newcommand{\ti}{\times}

\def\ga{\alpha}     \def\gb{\beta}       
       \def\gd{\delta}      \def\ge{\epsilon}
                         
\def\gf{\phi}           
            \def\gl{\lambda}
\def\gm{\mu}        \def\gn{\nu}         \def\gp{\pi}
    \def\gr{\rho}        
\def\gs{\sigma}       \def\gt{\tau}
      \def\gw{\omega}
                \def\gz{\zeta}
\def\Gg{\Gamma}     \def\Gd{\Delta}      
\def\Gth{\Theta}
          
\def\Gw{\Omega}              

      \def\CP{{\mathcal P}}

   \def\BBR {\mathbb R}

\allowdisplaybreaks
\begin{document}
\title {\bf  Propagation of Singularities of Nonlinear Heat Flow in Fissured Media}
\author{{\bf\large Andrey Shishkov}
\hspace{2mm}\vspace{3mm}\\
{\it \normalsize  Institute of Applied Mathematics and Mechanics of NAS of Ukraine},\\
{\it\normalsize R. Luxemburg str. 74, 83114 Donetsk, Ukraine}\\
\vspace{3mm}\\
{\bf\large Laurent V\'eron}
\vspace{3mm}\\
{\it\normalsize Laboratoire de Math\'ematiques et Physique Th\'eorique, CNRS UMR 6083},
\\
{\it\normalsize Universit\'e Fran\c{c}ois-Rabelais,  37200 Tours,
France}}

\date{}
\maketitle

{\it Dedicated to the memory of Igor Vladimirovich Skrypnik}
\section{Introduction}

In this paper we investigate the propagation of singularities in a nonlinear parabolic equation with strong absorption when the absorption potential is strongly degenerate following some curve in the $(x,t)$ space. As a very simplified model, we assume that the heat conduction is constant but the absorption of the media depends stronly of the characteristic of the media.
More precisely we suppose that the temperature $u$ is governed by the following equation
\begin{equation}\label{I-1}
\prt_{t}u-\Gd u+h(x,t)u^p=0\quad \text{in }Q_{T}:=\BBR^N\ti (0,T)
\end{equation}
where $p>1$ and $h\in C(\overline Q_{T})$. We suppose that $h(x,t)>0$ except when 
$(x,t)$ belongs to some space-time curve 
$\Gg$ given by 
\begin{equation}\label{I-2}
\Gg=\left\{\gamma (\gt):=(x(\gt),t(\gt)):\gt\in [0,T]\right\},
\end{equation}
where $\gamma\in C^{0,1}([0,1])$ with $\gamma(0)=0$ and $t(\gt)>0$ for $\gt\in (0,T]$. If there holds
\begin{equation}\label{I-3}
\myint{0}{T}\myint{B_{R}}{}h(x,t)E^q(x,t)dxdt<\infty
\end{equation}
we first show that, for any $k\geq 0$, there exists a unique solution $u_{k}$ of 
$(\ref{I-1})$ such that $u_{k}(.,0)=k\gd_{0}$. Furthermore the mapping $k\mapsto u_{k}$ is increasing. Because $h$ is continuous and positive outside $\Gamma$, we shall show that the set of solutions $\{u_{k}\}$ remains locally bounded in $Q_{T}\setminus \Gamma$. Therefore $u_{k}\uparrow u_{\infty}$ and $u_{\infty}$ is a solution of $(\ref{I-1})$ in $Q_{T}\setminus \Gamma$. Then either the singularity of the solution issued from the point $(0,0)$ can propagate along $\Gamma$ (at least partially), or it remains localized at $(0,0)$.
More precisely we show\medskip

\noindent  {\bf Theorem A.} {\it Assume $\gamma$ is $C^1$ and $t'(\gt)> 0$ for any $\gt\in (0,T]$. The following dichotomy of phenomena occurs

\noindent (i) either $u_{\infty}(x,t)<\infty$ for all $(x,t)\in Q_{T}$, \smallskip

\noindent (ii) or there exists $\gt_{0}\in (0,T]$ such that 
\begin{equation}\label{I-4}
\limsup_{(x,t)\to (y,s)}u_{\infty}(x,t)=\infty\quad \forall (y,s)\in\Gamma,0\leq s\leq \gt_{0}.
\end{equation}
}\medskip

We first prove that the singularity does not propagate from a point where $t(\gt)$ is decreasing. We prove\medskip

\noindent {\bf Theorem B.} {\it Assume that $\gamma$ is $C^1$ and $t'(\gt)<0$ for all $\gt\in (\gt_0,T]$ for some $\gt_0>0$. Then $u_\infty$ remains locally bounded in $Q_T\setminus \gamma((\gt_0,T])$.}\medskip

Due to this fact we shall assume first  that the $t$ variable is increasing along $\Gamma$, in such a case we can assume that $\gt$ is a function of $t$ and, up to a change of parameter, that $\gt=t$ and 
\begin{equation}\label{I-5}
\Gg=\left\{(x (t),t):t\in [0,T]\right\},
\end{equation}
where $\gamma\in C^{0,1}(0,T]$ satisfies $\gamma (0)=0$.
 In order the singularity to propagate along $\Gamma$, $h$ must be very flat near this curve $\Gamma$. If we define the parabolic distance between $(x,t)$ and $\Gamma$ by
$$d_{_{P}}(x,t);\Gamma)=\inf\{d_{_{P}}(x,t);(y,s)):(y,s)\in \Gamma,s\leq t\}$$ where $d_{_{P}}(x,t);(y,s)=|x-y|+\sqrt {t-s}$, for $t\geq s$
and write 
$$h(x,t)=e^{-\ell(d_{_{P}}(x,t);\Gamma))},
$$
where $\ell$ is a positive nonincreasing function. In the next result we shows that propagation of singularity along $\Gamma$ occurs

\smallskip{\it
\noindent {\bf Theorem C.} Assume $t\mapsto x(t)$ belongs to $W^{2,\infty}_{loc}[0,\infty)$ and
\begin{equation}\label{I-6}
\liminf_{t\to 0}t^2\ell(t)>0.
\end{equation}
Then $\lim_{(x,t)\to (y,s))}u_\infty(x,t)=\infty$ for all $(y,s)\in \Gamma$.
}\medskip

 The last section is devoted to the case where the curve of degeneracy $\Gamma$ is a straight line contained in the initial plan. We set $d_\infty((x,t),\Gamma)=\max\{\sqrt t,x'\}$ if $\dist (x,\Gamma)=x'$ and we write
 $$h(x,t)=e^{-\ell(d_{\infty}(x,t);\Gamma))}$$
 Up to a rotation, we can suppose that $\Gamma$ is the $x_1$ axis and denote $x=(x_1,x')$ the component in $\BBR\ti\BBR^{N-1}$. Then we prove the following\medskip
 
 \smallskip{\it
\noindent {\bf Theorem D.} Assume 
\begin{equation}\label{I-6}
\liminf_{t\to 0}t^2\ell(t)>0.
\end{equation}
Then $\lim_{(x,t)\to (y_1,0))}u_\infty(x,t)=\infty$ for all $y_1\in \BBR$.
}\medskip

Even if this model is a very simplified version of the heat propagation in a fissured absorbing media, it gives interesting insight of the propagation phenomenon which can occur. It is also a starting point for studying other type of propagation of singularities in nonlinear diffusion equations. In a forthcoming article we shall consider the case where the degeneracy line is a surface in $\BBR^N\ti [0,\infty)$ with only one contact point with $t=0$ at $(0,0)$.
\section{Preliminaries and basic estimates}
\setcounter{equation}{0}
Through out this section we assume that $h\in C(\BBR^N\ti[0,\infty))$ is nonnegative. 
Since $E(x,t)$ is a supersolution for $(\ref{I-1})$, the following result holds \cite[Theorem 6.12]{Ve1}
\bprop{Ex1} Assume $q>1$ and $(\ref{I-3})$ holds. Then for any $k>0$ there exists a unique $u=u_{k}\in C(\overline Q_{T}\setminus\{(0,0)\})\cap L^1(Q_{T})$,  such that $h|u|^q\in L^1(Q_{T})$, satisfying 
\begin{equation}\label{E-1}
\dint_{Q_{T}}(-u\prt_{t}\gz-\Gd\gz+h(x,t)|u|^{q-1}u\gz)dx\,dt
=k\gz (0,0)
\end{equation}
for all $\gz\in C^{2,1}(\overline Q_{T})$ which vanishes at $t=T$.
\es

More generaly, if $\gm\in \mathfrak M(\BBR^N)$ and $\gn\in \mathfrak M(\BBR^N\ti(0,T))$ are two positive bounded measures, the solution $v=v_{\gm,\gn}$ of
\begin{equation}\label{E-1-0}\left\{\BA {l}
\prt_{t}v-\Gd v=\gn\qquad\text{in } \BBR^N\ti[0,T)\\
\phantom{-\Gd }
v(.,0)=\gm\qquad\text{in } \BBR^N,
\EA\right.\end{equation}
is expressed by
\begin{equation}\label{E-1-1}
v_{\gm,\gn}(x,t)=\myfrac{1}{(4\gp t)^{N/2}}
\myint{\BBR^N}{}e^{-4|x-y|^2/4t}d\gm(y)+
\myint{0}{t}\myfrac{1}{(4\gp (t-s))^{N/2}}
\myint{\BBR^N}{}e^{-4|x-y|^2/4(t-s)}d\gn(y)ds.
\end{equation}
Actually, by direct adaptation to the parabolic case of \cite[Theorem 4.2]{Ve1}, combined with \cite[Theorem 6.12]{Ve1},  one can proves that, for any bounded  Radon  measures $\gm$ on $\BBR^N\ti (0,\infty)$ and $\gn$ on $\BBR^N$ which satisfy
\begin{equation}\label{E-1-3}
\dint_{Q_{T}}v^q_{\gm_{s},\gn_{s}}h(x,t)dxdt<\infty
\end{equation}
where $\gm_{s}$ and $\gn_{s}$ are the singular parts (with respect to respective Lebesgue measures) $\gm$ and $\gn$ respectively, there exists a unique weak solution $u$ to problem
\begin{equation}\label{E-1-2}\left\{\BA {l}
\prt_{t}u-\Gd u+h(t,x)|u|^{q-1}u=\gn\qquad\text{in } \BBR^N\ti[0,T)\\
\phantom{-\Gd+h(t,x)|u|^{q-1}u }
u(.,0)=\gm\qquad\text{in } \BBR^N,
\EA\right.\end{equation}
The next result is an adaptation of Brezis-Friedman a priori estimate \cite{BF}.
 \bprop{Ex2} Let $Q=Q^{r,a}_{t_{0},t_{1}}:= B_{r}(a)\ti (t_{0},t_{1})$ for some 
 $a\in\BBR^N$, $t_{1}>t_0\geq 0$ and $r>0$ and assume 
 $\beta=\min\{h(x,t):(x,t)\in \overline Q \}>0$.
Then any solution of $(\ref{I-1})$ in $Q$ satisfies
 \begin{equation}\label{E-2}
|u(x,t)|\leq \myfrac{C}{\beta^{1/(q-1)}}
\left(\myfrac{1}{t-t_{0}}
+\myfrac{1}{(r-|x-a|)^2}\right)^{1/(q-1)}\qquad\forall (x,t)\in Q,
\end{equation}
for some $C=C(N,q)>0$.
 \es
 \Proof The maximal solution of the parabolic equation
 $y'+\beta y^q=0$ on $0,\infty)$ is expressed by 
 \begin{equation}\label{y_M}y_{_{M}}(t)=\left(\myfrac{1}{\beta (q-1)t}\right)^{1/(q-1)}.\end{equation}
 By Keller-Osserman estimate, the maximal solution $v_{_{M}}$ of $-\Gd v+\beta |u|^{q-1}u=0$ in $B_{r}$, satisfies
 \begin{equation}\label{v_M}v_{_{M}}(x)\leq C\left(\myfrac{1}{\beta (r-|x|)^2}\right)^{1/(q-1)}
\end{equation}
 in $B_{r}$. Since $y_{_{M}}(t-t_{0})+v_{_{M}}(x-a)$ is a supersolution of $(\ref{I-1})$ in $Q$ which blows up on the parabolic boundary, an easy approxition argument (just replacing $r$ by $\{r_{n}\}\uparrow r$ and $t_{0}$ by $\{t_{n}\}\downarrow t_{0}$ ) leads us to $(\ref{E-2})$.\qeda\medskip

 \bprop{Ex3} Let $0<r_0<r_1$, $0<t_0<t_1$ and $\Theta:=\Theta^{r_0,r_1,a}_{t_0,t_1}:=Q^{r_1,a}_{t_{0},t_{1}}\setminus Q^{r_0,a}_{t_{0},t_{1}}$ for some $a\in\BBR^N$ and assume $\beta=\min\{h(x,t):(x,t)\in \overline \Theta \}>0$. Then any solution of (\ref{I-1}) in $Q^{r_1,a}_{t_{0},t_{1}}$ such that $|u(x,t_0)|\leq \gm$ for $x\in B_{r_1}(a)$ satisfies
 \begin{equation}\label{E-2-1}
|u(x,t_1)|\leq \gm+\myfrac{C}{(\beta(r_1-r_0)^2)^{1/(q-1)}}\qquad\forall x\in B_{r_0}(a),
\end{equation}
 for some $C=C(N,q)>0$.
\es
\Proof Let $b\in\BBR^N$ such that $|b-a|=\frac{r_0+r_1}{2}$ and $v_M$ the maximal solution of 
$-\Gd v+\beta |u|^{q-1}u=0$ in $B_{\frac{r_1-r_0}{2}}(b)$. Then
$$
v_M(x)\leq \myfrac{C}{\left(\beta\left(\frac{r_1-r_0}{2}-|x-b|\right)\right)^{2/(q-1)}}\qquad\forall x\in B_{\frac{r_1-r_0}{2}}(b).
$$
Since $v_M+\gm$ is a super solution of (\ref{I-1}) in $Q^{\frac{r_1-r_0}{2},b}_{t_{0},t_{1}}$ which dominates $u$ at $t=t_0$ and for $|x-b|=\frac{r_1-r_0}{2}$, it follows that $u(x,t)\leq v_M(x,t)+\gm$ in $Q^{\frac{r_1-r_0}{2},b}_{t_{0},t_{1}}$. 
In particular, if $x=b$, we get
$$
u(b,t)\leq \gm+\myfrac{C}{\left(\beta\left(\frac{r_1-r_0}{2}\right)\right)^{2/(q-1)}}
$$
This estimate is valid for any $b\in\BBR^N$ with $|b-a|=\frac{r_0+r_1}{2}$. Since $u$ is a subsolution of the heat equation in 
$Q^{r_0,a}_{t_0,t_1}$, $(\ref{E-2-1})$ follows by the maximum principle.
\qeda\medskip

The previous estimates are based upon constructions of supersolutions in cylinders. In the next result we construct estimates in tubular neighborhood of $\Gamma$. If $t(\gt)$ is increasing, we can take $\gt=t$ and 
\begin{equation}\label{incr}
\Gamma=\gamma ([0,T]):=\{(x(t),t):t\in [0,T]\}
\end{equation}

\bprop{Ex4} Assume $\Gamma$ is $C^1$ and parametrized by $t$ as in $(\ref{incr})$, and for $\ge>0$ denote 
$$ T^{\Gg_0^T}_s:=\{(x,t): 0\leq t\leq T,\,|x-x(t)|<s\}.
$$
For $1\geq r_1>r_0>0$ we set $\eta=\min\left\{h(x,t):(x,t)\in T^{\Gg_0^T}_{r_1}\setminus T^{\Gg_0^T}_{r_0}\right\}$. Then there exists a constant $C$ depending on $c:=\max \{|x'(t)|\ :0\leq t\leq T\}$ such that if 
\begin{equation}\label{incr'}
u(x,0)\leq m+\frac{C_2}{(r_1-r_0)^{\frac{2}{q-1}}}
\qquad\forall x\in \overline{B_{\frac{r_1+r_0}{2}}(x(0))},
\end{equation}
then
\begin{equation}\label{incr0}
u(x,t)\leq \left(\frac{m^{q-1}}{1+\eta(q-1)m^{q-1}t}\right)^{\frac{1}{q-1}}+\frac{C_2}{(r_1-r_0)^{\frac{2}{q-1}}}
\qquad\forall t\in [0,T] \text { and }x\in \overline{B_{\frac{r_1+r_0}{2}}(x(t))}
\end{equation}
\es
\Proof Consider the change of space variable $x=y+x(t)$ and $u(x,t)=v(y,t)$. Then $v$ satisfies
\begin{equation}\label{incr1}
\prt_tv-\Gd v+\langle x'(t),\nabla v\rangle +h(y+x(t),t)v^q=0\qquad\text{ in }Q_{0,t(T)}^{r_0,r_1}.
\end{equation}
Thus 
$$\prt_tv-\Gd v-c|\nabla v| +\eta v^q\leq 0\qquad\text{ in }Q_{0,t(T)}^{r_0,r_1}.
$$

In the ball $B_{\frac{r_1-r_0}{2}}(z)$ where $|z|=\frac{r_1+r_0}{2}$ It is standard easy to construct a radial function $\psi$ satisfying
\begin{equation}\label{incr2}-\Gd \psi-c|\nabla \psi| +\eta \psi^q\geq 0 \qquad\text{in }B_{\frac{r_1-r_0}{2}}(z)
\end{equation}
under the form
$$\psi (y)=C\frac{\gr^{\frac{2}{q-1}}}{(\gr^2-|y-z|^2)^{\frac{2}{q-1}}}
$$
where $\gr=\frac{r_1-r_0}{2}$ and $C=C(N,q,c,\eta)$. Therefore
$$\psi(z)=\frac{C_1}{\gr^{\frac{2}{q-1}}}=\frac{C_2}{(r_1-r_0)^{\frac{2}{q-1}}}
$$
The solution $\gf_m=\gf$ of
\begin{equation}\label{incr3}\left\{\BA {l}
\gf'+\eta\gf^q=0\qquad\text {on } (0,T)\\
\phantom{\gf'\gf^q}
\gf(0)=m>0
\EA\right.\end{equation}
is expressed by
$$\gf_m(t)=\left(\frac{m^{q-1}}{1+\eta(q-1)m^{q-1}t}\right)^{\frac{1}{q-1}}.
$$
Consequently, the function $\gf_m(t)+\psi (y)$ is a supersolution for $(\ref{incr1})$ in  $Q_{0,t(T)}^{r_0,r_1}$ which domintates $v$ at $t=0$ and on the lateral boundary. Therefore it is larger than $v$ and in particular
\begin{equation}\label{incr4}
v(z,t)\leq \gf_m(t)+\psi(z).
\end{equation}
Consequently
\begin{equation}\label{incr4}
u(x,t)\leq \gf_m(t)+\psi(z)\qquad\forall t\in (0,T) \text { and } |x-x(t)|=\frac{r_1+r_0}{2}.
\end{equation}
We derive $(\ref{incr0})$ by the maximum principle.\qeda

\section{Geometric obstruction to propagation}
\setcounter{equation}{0}
\medskip

We assume that $h$ vanishes on a continuous curve $\Gamma\subset \BBR^N\ti [0,T)$ defined by parametrisation
 \begin{equation}\label{E-3}
\Gamma=\{\gamma(t):=(x(\gt),t(\gt):\gt\in [0,T]\}
\end{equation}
issued from $(0,0)$ (i.e. $x(0),t(0)=(0,0)$ with $t(\gt)>0$ if $\gt\in (0,T)$. and $\gt\mapsto\gamma (\gt)$ is Lipschitz with no self intersection, which means that $\gt\mapsto\gamma(\gt)$ is one to one. \medskip


\noindent{\it Proof of Theorem A.} Since $t'$ is continuous and increasing, we apply \rprop {Ex4} with $0$ replaced by $\gt_0$, we set $t_0=t(\gt_0)$ and write $\Gamma$ under the form $(\ref{incr})$. If we assume that
$$\limsup_{(x,t)\to(x(t_0),t_0)}u(x,t)<\infty,
$$
there exists $r_1>0$ and $\gm>0$ such that $u(x,t_0)\leq \gm$ if $x\in B_{r_1}(x(t_0))$. Then, for $0<r_0<r_1$ there exists 
$m>0$ such that $(\ref{incr})$ is verified. Thus $(\ref{incr1})$ holds. This implies that the blow-up set of $u$ along $\Gamma$ is empty if $t\geq t_0$. 
\qeda\medskip


\noindent{\it Proof of Theorem B. }The proof is based upon the same ideas than in Theorem A above except that we only study the part of $\Gamma$ between $\gt_0$ and $T$, where $t'(\gt)<0$. Let $t_0=t(\gt_0)$ and  $t^*=t(T)$. We parametrized $\Gamma$ by $t$ between $t^*$ and $t_0$, thus 
$$\gamma([\gt_0,T])=\{(t,x(t)):t^*\leq t\leq t_0\}.
$$
For $s>0$ and $t^*< t'\leq t_0$, we set $T_s^{\Gamma_{t*}^{t'}}:=\{(x,t): t^*\leq t\leq t_0,\,|x-x(t)|< s\}$, then, for  $t^*< t_1< t_0$, there exists $r_1>0$ and $\gt_1\in (\gt_0,T)$ such that $t_1=t(\gt_1)$ and $\Gamma\cap T_{r_1}^{\Gamma_{t*}^{t_1}}=\gamma ([\gt_1,T])$. If $t^*=0$, then $u(x,t^*)=0$ for $|x-x(t^*)|\leq r_1$. If $t^*>0$, then $h(x,t)>0$ in the cylinder $Q^{r_1,x(t^*)}_{0,t^*}$, thus there exists $\beta>0$ such that $\inf\left\{h(x,t):(x,t)\in Q^{r_1,x(t^*)}_{0,\frac{t^*}{2}}\right\}=\beta$. Up to replacing $r_1$ by some $r'_1>r_1$ such that $\Gamma\cap \overline {Q^{r_1,x(t^*)}_{0,\frac{t^*}{2}}}=\emptyset$, we obtain that $u_\infty(x,t^*)\leq \gm$ for all $x\in \overline {B_{r_1}(x(t^*))}$ from \rprop{Ex2}. In both cases $u_\infty(x,t^*)$ is bounded in $\overline {B_{r_1}(x(t^*))}$. Replacing $(0,T)$ by $(t^*,t_1)$) in \rprop{Ex4}, it follows that $u_\infty$ remains bounded in $T_{\frac{r_1+r_0}{2}}^{\Gamma_{t*}^{t_1}}$ for some $0<r_0<r_1$. Since $t_1<t_0$ is arbitrary and $u_\infty$ is locally bounded in $Q_T\setminus \Gamma$, the proof follows.
\qeda\medskip

In the next case the monotonicity of $\gt\mapsto t(\gt)$ is replaced by a {\it box-assumption}.
\bprop {GBox}Assume that $\gamma$ is continuous and there exists $a\in\BBR^N$, $r_0>0$ and $\gt_0\in (0,T)$ such that $t(T)\leq t(\gt_0)$ and $\gamma([\gt_0, T])\subset B_{r_0}(a)\ti [t(T),t(\gt_0)]$. Then 
$u_{\infty}$ is bounded in $\overline B_{r_0}(a)\ti [t(T),t(\gt_0)]$.
\es
\Proof  There exist $r'_0<r_0$ and $\beta>0$ such that $\gamma([\gt_0, T])\subset B_{r'_0}(a)\ti [t(T),t(\gt_0)]$ and 
$\min\left\{h(x,t):(x,t)\in \Gth^{r'_0,r_0}_{t(T),t(\gt_0)}\right\}=\beta$. Therefore the conclusion follows from \rprop{Ex3}.\qeda\medskip

In the next case we show that non-propagation of singularities may occur even if $\gt\mapsto t(\gt)$ is increasing after some $\gt_0$, provided there is a local maximum in $(0,\gt_0)$. We put $\Gg_{\gt_0}=\gamma([0,\gt_0])$

\bth {GBoxes}Let $\gamma$ be $C^1$ and $\gamma'(\gt)\neq 0$ on $[0,T]$. Assume there exist $\gt_0>0$, $a\in\BBR^N$ and $r>0$ such that $t(\gt)\leq t(\gt_0)$ on $[0,\gt_0]$, $\gt\mapsto t(\gt)$ is decreasing on $(\gt_0,\gt_0+\gd)$ for some $\gd>0$, $\gamma((\gt_0,T])\subset Q_{r,a}^{0,\infty}$,  $\Gg_{\gt_0}\cap \overline{Q_{r,a}^{0,\infty}}=\{\gamma(\gt_0)\}$ and  for any $\gt \in (\gt_0,T]$, $\abs{x(\gt)-a}<r$. Then $u_{\infty}$ is locally bounded in $Q\setminus\Gg_{\gt_0}$, where $\Gg_{\gt_0}:=\gamma([0,\gt_0])$.
\es
\Proof Since $t(\gt)$ is decreasing on $(\gt_0,\gt_0+\gd)$, for any $\gt_1\in (\gt_0,\gt_0+\gd)$, the set of $\gt\in (\gt_1,T]$ such that $t(\gs)<t(\gt_1)$ for all $\gs\in (\gt_1,\gt)$ is not empty. Its upper bound $\gt^*_1$ is less or equal to $T$ and $t(\gt^*_1)=\min\{t(\gt_1),t(T)\}$. If $\lim_{\gt_1\to\gt_0}t(\gt^*_1)=t(T)$, then $t(\gt_0)=t(T)$, the {\it box-assumption} holds and the conclusion follows from \rprop{GBox}. If $t(T)>t(\gt_0)$, then $t(\gt_1)=t(\gt^*_1)<t(T)$ for any $\gt_1\in (\gt_0,\gt_0+\gd)$. Since 
$\lim_{\gt_1\to\gt_0}t(\gt_1)=t(\gt_0)$ and $\gamma$ is continuous with $\gamma((\gt_0,T])\subset Q_{r,a}^{0,\infty}$, there some fixed constants $\lambda>0$ and $\gr>0$ such that $x(\gt)\in B_{\gr}(x(\gt_1))$, for any $\gt\in [\gt_1,T]$ verifying $t(\gt)\leq t(\gt_1)+\lambda$. This means that the part of $\Gg$ starting from $\gamma (\gt_1)$ for which $t(\gt)$ belongs to $(t(\gt_1), t(\gt_1)+\lambda)$ remains in $B_{\gr}(x(\gt_1))$. Moreover we can assume that $t(\gt_1)+\lambda>t(\gt_0)$. By restricting $\gr$, we can assume that $B_{\gr}(x(\gt_1))\subset B_{r'}(a)$. Since $u_k(x,t(\gt_1)$ is uniformly bounded when $x\in B_{\gr}(x(\gt_1))$, it follows from 
\rprop{Ex3} that for any $\gr'\in (0,\gr)$, $u_k$ remains uniformly bounded in $Q^{t(\gt_1),t(\gt_1)+\lambda}_{\gr',x(\gt_1)}$. Moreover, for any compact subset $K\in B^c_{\gr'}(x_{\gt_1})$, there holds
$$u_k(x,t(\gt_1)+\lambda)\leq \myfrac{C}{(\lambda\beta)^{1/(q-1)}},
$$
where $\gb=\min\{h(x,t):(x,t)\in K\ti [t(\gt_1),t(\gt_1)+\lambda]\}$.
Iterating this construction, we can construct a finite number of cylinders $Q^{t(\gt_1)+j\lambda,t(\gt_1)+(j+1)\lambda}_{\gr',a_j}$ containing  $\Gamma$ and in which $u_k$ remains uniformly bounded. Since local uniform boundedness holds also outside such  cylinders, the proof follows.
\qeda\medskip

\noindent\Remark In full generality we conjecture that  if $\gamma$ is $C^1$ with $\gamma'(\gt)\neq 0$ and there exists $\gt_{0}\in (0,T)$ such that $t(\gt)$ admits a local strict maximum on the right at  $\gt_{0} \in  (0, T)$, then
$u_{\infty}$ is locally bounded $ Q_{T}\setminus \Gamma_{\gt_{0}}$.

\medskip
\section{Propagation of singularities in the space}
\setcounter{equation}{0}
In this section, we assume that the degeneracy curve $\Gamma$ is parametrized by the variable $t\in [0,T]$  and defined by $(\ref{incr})$ with $|x (t)|>0$ if $t>0$. 
We denote by 
$$d_{_{P}}[(x,t),(y,s)]:=|x-y|+\sqrt{t-s}\quad \text{if }t\geq s$$
the parabolic distance and we assume that $h(x,t)$ depends on $d_{_{P}}[(x,t),\Gg]$ under the following form
\begin{equation}\label{cur2}
h(x,t)=e^{-\ell (d_{_{P}}[(x,t),\Gg])}
\end{equation}
where $\ell\in C([0,\infty))$ is positive, nonincreasing and $\lim_{r\to 0}\ell (r)=\infty$. For $\ge>0$, we recall 
$T^{\Gg_0^T}_{\ge}$ denotes the $\ge$ -spherical tubular neighborhood of $\Gg$ between $t=0$ and $t=T$ defined by 
\begin{equation}\label{cur3}
T^{\Gg_0^T}_{\ge}:=\{(x,t)\in Q_{T}:|x-x(t)|<\ge\}.
\end{equation}
The basis of $T^{\Gg_0^T}_{\ge}$, in $\BBR^N\ti\{0\}$, is the ball $B_{\ge}$. 
Since $d_{_{P}}[(x,t),\Gg]\leq |x-x(t)|$, $d_{_{P}}[(x,t),\Gg]\leq \ge$ in $T^{\Gg_0^T}_{\ge}$ and 
$\ell(\ge)\leq \ell(d_{_{P}}[(x,t),\Gg])$. Then, $u_\infty$ is bounded from below in 
$\Gg_{\ge}$ by the solution $v_{\ge}$ of
\begin{equation}\label{cur4}
\left\{\BA {l}
\prt_{t} v_{\ge}-\Gd v_{\ge}+e^{-\ell(\ge)}v^p_{\ge}=0\quad\quad\text {in }T^{\Gg_0^T}_{\ge}\\[2mm]
\phantom{\prt_{t} v_{\ge}-\Gd v_{\ge}^p+e^{-\ell(\ge)}}
v_{\ge}=0\quad\quad\text {in }\prt_{_{\CP}}T^{\Gg_0^T}_{\ge}\\[2mm]
\phantom{\prt v_{\ge}-\Gd v_{\ge}^p+e^{-\ell(\ge)}}
v_{\ge}=\infty\gd_{0}\quad\text {in }B_{\ge},
\EA\right.
\end{equation}
where 
$$\prt_{_{\CP}}T^{\Gg_0^T}_{\ge}:=\{(x,t):|x-x(t)|=\ge\}$$ 
is the lateral parabolic boundary of $T^{\Gg_0^T}_{\ge}$. The formal aspect comes from the fact that the existence of $v_{\ge}$ has to be proved. \smallskip

{\it For the sake of simplicity we assume that $1<p<1+\frac{2}{N}$, the case $p\geq 1+\frac{2}{N}$ needing a simple adaptation of the type which is developed in the proof  of Theorem D-Case 2}. \smallskip 

We consider the following change of variable $x=y+x(t)$ and
$v_{\ge}(x,t)=\tilde v_{\ge}(y,t)$. Then $\tilde v_{\ge}$ satisfies in $Q_{\gt}^{B_{\ge}}:=B_{\ge}\ti (0,\gt)$
\begin{equation}\label{cur5*}
\left\{\BA {l}
\prt_{t} \tilde v_{\ge}-\Gd \tilde v_{\ge}+\langle x'(t)|\nabla \tilde v_{\ge}\rangle+e^{-\ell(\ge)}\tilde v^p_{\ge}=0\quad\quad\text {in }Q_{\gt}^{B_{\ge}}\\[2mm]
\phantom{\prt_{t} \tilde v_{\ge}-\Gd +\langle x'(t)|\nabla \tilde v_{\ge}\rangle+e^{-\ell(\ge)}\tilde v^p_{\ge}}
 \tilde v_{\ge}=0\quad\quad\text {in }\prt_{_{\CP}}Q_{\gt}^{B_{\ge}}\\[2mm]
\phantom{\prt_{t} \tilde v_{\ge}-\Gd +\langle x'(t)|\nabla \tilde v_{\ge}\rangle+e^{-\ell(\ge)}\tilde v^p_{\ge}}
 \tilde v_{\ge}=\infty\gd_{0}\quad\text {in }B_{\ge}.
\EA\right.
\end{equation}
 In particular 
$u(x (\ga),\ga)\geq v_{\ge}(x (\ga),\ga)=\tilde v_{\ge}(0,\ga)$. We set
$$\gw_{\ge}(x,t)=\ge^{2/(p-1)}e^{-\ell(\ge)/(p-1)} \tilde v_{\ge}
(\ge x,\ge^2t),
$$
and $x_{\ge}(t)=\ge^{-1}x(\ge^2t)$.
Then $\gw_{\ge}$ satisfies 
\begin{equation}\label{cur5**}
\left\{\BA {l}
\prt_{t} \gw_{\ge}-\Gd \gw_{\ge}+\langle x_{\ge}'(t)|\nabla \gw_{\ge}\rangle+\gw^p_{\ge}=0\quad\quad\text {in }Q_{\ge^{-2}\gt}^{B_{1}}\\[2mm]
\phantom{\prt_{t} \gw_{\ge}-\Gd+\langle x_{\ge}'(t)|\nabla \gw_{\ge}\rangle+\gw^p_{\ge}}
\gw_{\ge}=0\quad\quad\text {in }\prt_{_{\CP}}Q_{\ge^{-2}\gt}^{B_{1}}\\[2mm]
\phantom{\prt_{t} \gw_{\ge}-\Gd+\langle x_{\ge}'(t)|\nabla \gw_{\ge}\rangle+\gw^p_{\ge}}
\gw_{\ge}=\infty\gd_{0}\quad\text {in }B_{1},
\EA\right.
\end{equation}
and for any $\ga>0$
\begin{equation}\label{cur5'}u_\infty(x (\ga),\ga)\geq \ge^{-2/(p-1)}e^{\ell(\ge)/(p-1)}\gw_{\ge}(0,\ge^{-2}\ga).
\end{equation}
Therefore, the problem is reduced to showing that
\begin{equation}\label{cur6}
\lim_{\ge\to 0}\ge^{-2/(p-1)}e^{\ell(\ge)/(p-1)}\gw_{\ge}(0,\ge^{-2}\ga)=\infty.
\end{equation}
We associate the following parabolic equation
\begin{equation}\label{cur5}
\left\{\BA {l}
\prt_{t} v-\Gd v+\langle\gb(t)|\nabla v\rangle+v^p=0\quad\quad\text {in }Q_{\infty}^{B_{1}}\\[2mm]
\phantom{\prt_{t} v-\Gd +\langle\gb(t)|\nabla v\rangle+v^p}
v=0\quad\quad\text {in }\prt_{_{\CP}}Q_{\infty}^{B_{1}}\\[2mm]
\phantom{\prt_{t} v-\Gd +\langle\gb(t)|\nabla v\rangle+v^p}
v(y,0)=v_{0}\quad\text {in }B_{1}.
\EA\right.
\end{equation}
Existence of solution is classical when $\gb\in L_{loc}^{\infty}([0,\infty))$ (see \cite{Fr}, \cite{LSU}).
\bprop{asympt1} Let $\gb\in L_{loc}^{\infty}([0,\infty))$. Then the following estimates holds
\begin{equation}\label{est1}
\norm{v(.,t)}_{L^2}\leq e^{-\gl_{0}t}\norm{v(.,0)}_{L^2}
\end{equation}
where $\gl_{0}$ is the first eigenvalue of $-\Gd$ in $B_{1}$, and, for some $C=C(N)>0$,
\begin{equation}\label{est2}
\norm{v(.,t)}_{L^\infty}\leq C\min\{t^{-N/4},e^{-\gl_{0}t}\}\norm{v(.,0)}_{L^2}.
\end{equation}
\es
\Proof
Since for any continuous function $r\mapsto g(r)=G'(r)$, we have
$$\myint{B_{1}}{}\langle\gb(t)|\nabla v\rangle g(v)dx=\myint{B_{1}}{}\langle\gb(t)|\nabla G(v)\rangle dx=-\myint{B_{1}}{} G(v)\,div\gb(t)dx=0,
$$
we have
$$2^{-1}\myfrac{d}{dt}\myint{B_{1}}{}v^2dx+\myint{B_{1}}{}|\nabla v|^2dx\leq 0\Longrightarrow \myfrac{d}{dt}\myint{B_{1}}{}v^2dx\leq -2\gl_{0}
\myint{B_{1}}{}v^2dx.$$
Thus $(\ref{est1})$ follows. Furthermore, for any $q>1$,
$$\myfrac{1}{q+1}\myfrac{d}{dt}\myint{B_{1}}{}|v|^{q+1}dx+\myfrac{4q}{(q+1)^2}
\myint{B_{1}}{}|\nabla v^{(q+1)/2}|dx\leq 0,
$$
it follows that $t\mapsto \norm{v(.,t)}_{L^{q+1}}$ is decaying. Therefore, using Gagliardo-Nirenberg estimate, 
$$\myfrac{d}{dt}\norm{v}^{q+1}_{L^{q+1}}+\myfrac{4qC}{q+1}
\norm{v}^{q+1}_{L^{(q+1)N/(N-2)}}\leq 0.
$$
We assume $N\geq 3$, the cases $N=1,2$ needing a simple modification. Thus, for any $t>s>0$,
$$\myfrac{4qC(t-s)}{q+1}
\norm{v(.,t)}^{q+1}_{L^{(q+1)N/(N-2)}}\leq \norm{v(.,s)}^{q+1}_{L^{q+1}}.
$$
By a standard use of Moser iterative scheeme, we derive
\begin{equation}\label{est3}
\norm{v(.,t)}_{L^\infty}\leq Ct^{-N/4}\norm{v(.,0)}_{L^2}.
\end{equation}
Consequently, for any $0<s<t$, we have
$$\BA {l}
\norm{v(.,t)}_{L^\infty}\leq C(t-s)^{-N/4}\norm{v(.,s)}_{L^2}\\[2mm]
\phantom{\norm{v(.,t)}_{L^\infty}}
\leq 
C(t-s)^{-N/4}e^{-\gl_{0}s}\norm{v(.,0)}_{L^2}.
\EA$$
Noting that, if $t>N/4\gl_{0}$, 
$$\min\{(t-s)^{-N/4}e^{-\gl_{0}s}0<s<t\}=e^{N/4}(4\gl_{0}/N)^{N/4}e^{-\gl_{0}t},$$
we derive $(\ref{est2})$.
\qeda\medskip

\bprop{asympt2} Let $\Gw$ be a bounded open domain in $\BBR^N$, $\gb\in\BBR^N$ and $L_{\gb}$ the operator $v\mapsto -\Gd v+\langle \gb|\nabla v\rangle$. Then the spectrum of $L_{\gb}$ in $H^{1,0}(\Gw)$ is the given by
\begin{equation}\label{spec}
\gs(L_{\gb})=\left\{\gl+\frac{|\gb|^2}{4}:\gl\in \gs(L_{0})\right\}.
\end{equation}
\es
\Proof  Put $v(x)=e^{\frac{1}{2}\langle \gb|x\rangle}w(x)$. Then
$$\nabla v(x)=e^{\frac{1}{2}\langle \gb|x\rangle}\left(\nabla w+\frac{w}{2}\gb\right)
$$
$$\Gd v(x)=e^{\frac{1}{2}\langle \gb |x\rangle}\left(\Gd w+\langle \gb |\nabla w\rangle+\frac{|\gb|^2}{4}w\right)
$$
Thus
\begin{equation}\label{spec1}
L_{\gb}v=e^{\frac{1}{2}\langle \gb |x\rangle}\left(-\Gd w+\frac{|\gb|^2}{4}w\right),
\end{equation}
and the proof follows.
 \qeda\medskip

In the sequel we denote by  $\gl_{\gb}$ the first eigenvalue of $L_{\gb}$, thus $\gl_\gb=\gl_0+\frac{|\gb|^2}{4}$. If $\psi_\gb$ is a corresponding positive eigenfunction, then 
$$\psi_\gb(x)=e^{\frac{1}{2}\langle \gb|x\rangle}\psi_0(x)
$$
where $\psi_0$ is a positive first eigenfunction of $-\Gd$ in $H^{1,0}(\Gw)$

\bprop{aut} Under the assumptions of \rprop{asympt2} and for $p>1$, we denote by $v$ be the solution of 
\begin{equation}\label{cur7}
\left\{\BA {ll}
\prt_{t} v-\Gd v+\langle\gb|\nabla v\rangle+v^p=0\quad\quad&\text {in }Q_{\infty}^{B_{1}}\\[2mm]
\phantom{\prt_{t} v-\Gd +\langle\gb|\nabla v\rangle+v^p}
v=0\quad\quad&\text {in }\prt_{_{\CP}}Q_{\infty}^{B_{1}}\\[2mm]
\phantom{\prt_{t} ,, +\langle\gb|\nabla v\rangle+v^p}
v(y,0)=v_{0}\quad&\text {in }B_{1}.
\EA\right.
\end{equation}
where $v_{0}\in L^2(B_{1})$ is nonnegative. Then there exists some $c=c(v_{0})>0$ such that
\begin{equation}\label{lim1}
\lim_{t\to\infty}e^{\gl_{\gb}t}v(.,t)=c\psi
\end{equation}
uniformly in $B_{1}$.
\es
\Proof We write $v(x,t)=e^{\frac{1}{2}\langle \gb|x\rangle}w(x,t)$, thus $(\ref{cur7})$ turns into
\begin{equation}\label{lim2}
\prt_{t} w-\Gd w+\frac{|\gb|^2}{4}w+e^{\frac{p-1}{2}\langle \gb|x\rangle}w^p=0.
\end{equation}
The proof of \cite[Th 3.1]{GV} applies easily and the result follows.\qeda

\bprop{asympt3} Assume $\gb\in W^{1,\infty}_{loc}([0,\infty))$. For $\gt>1$, we set
$\sup_{1\leq t\leq \gt}|\gb(t)|=\gb_\gt$ and $\sup_{1\leq t\leq \gt}|\gb'(t)|=\gd_\gt$. If $v$ is the solution of $(\ref{cur5})$ where $v_{0}\in L^2(B_{1})$ is nonnegative, we denote
 $\gs(\gt)=\sup\left\{e^{\frac{p-1}{2}|\gb(t)|}w^{p-1}(x,t):(x,t)\in B_1\ti [1,\gt]\right\}$. Then if $v(.,1)$ satisfies $c_1\psi_0\leq v(.,1)$, there holds
\begin{equation}\label{W0}
 v(x,t)\geq c_1 e^{-(\gl_0+\frac{\gb_\gt^2}{4}+\frac{\gd_\gt}{2}+\gs_\gt)(t-1)}\psi_{0}(x)\quad\forall (x,t)\in B_1\ti [1,\gt]
\end{equation}
\es
\Proof Let $w(x,t)=e^{-\frac{1}{2}\langle \gb(t)|x\rangle}v(x,t)$. Then 
\begin{equation}\label{W1}
\prt_{t} w-\Gd w+\left(\frac{|\gb|^2}{4}+\frac{1}{2}\langle\gb'|x\rangle+e^{\frac{p-1}{2}\langle \gb(t)|x\rangle}w^{p-1}\right)w=0.
\end{equation}
Since $|x|\leq 1$ there holds in $B_1\ti [1,\gt]$
\begin{equation}\label{W2}
\frac{|\gb(t)|^2}{4}+\frac{1}{2}\langle\gb'(t)|x\rangle+e^{\frac{p-1}{2}\langle \gb(t)|x\rangle}w^{p-1}\leq \frac{\gb_\gt^2}{4}+\frac{\gd_\gt}{2}+\gs_\gt.
\end{equation}
Therefore
\begin{equation}\label{W2}
\prt_{t} w-\Gd w+\left(\frac{\gb_\gt^2}{4}+\frac{\gd_\gt}{2}+\gs_\gt\right)w\geq 0\quad\text{in }B_1\ti [1,\gt].
\end{equation}
Since $(x,t)\mapsto e^{-(\gl_0+\frac{\gb_\gt^2}{4}+\frac{\gd_\gt}{2}+\gs_\gt)(t-1)}\psi_{0}(x)$ satisfies the equation associated to $(\ref{W2})$, 
$(\ref{W0})$ follows.\qeda\medskip

\smallskip

\noindent{\it Proof of Theorem C.} {\it Step 1: Initialization of the blow-up.} With our previous notations, $\gb(t)=\gb_\ge(t)=x'_\ge(t)=\ge x'(\ge^2t)$ and $\gb_\ge'(t)=\ge^3 x''(\ge^2t)$. Since $x'_\ge$ is locally bounded, it follows by Hopf lemma that there exists $c_1>0$ such that 
$\gw_\ge(.,1)\geq c_1\psi_0$. We take $\gt=\ge^{-2}\ga$ where $\ga>0$ is fixed. Then
$$\gb_\gt=\sup\{|\ge x'(\ge^2t)|:1\leq t\leq \ge^{-2}\ga\}=\ge\sup\{| x'(t)|:\ge^2\leq t\leq \ga\},$$
$$\gd_\gt=\sup\{|\ge^3 x''(\ge^2t)|:1\leq t\leq \ge^{-2}\ga\}=\ge^3\sup\{| x''(t)|:\ge^2\leq t\leq \ga\}.
$$
Therefore
\begin{equation}\label{W3}\ge^{-\frac{2}{p-1}}e^{\frac{\ell(\ge)}{p-1}}\gw_\ge(0,\ge^{-2}\ga)\geq c'_1e^{A(\ge)}\psi_0(0)
\end{equation}
where
\begin{equation}\label{W4}A(\ge):=-\frac{2}{p-1}\ln\ge+\frac{\ell(\ge)}{p-1}-\left(\gl_0+\frac{\gb_\gt^2}{4}+\frac{\gd_\gt}{2}+\gs_\gt\right)\frac{\ga}{\ge^{2}}
\end{equation}
Since  $\liminf_{\ge\to 0}\ge^2\ell (\ge)>0$
there exists $\ga_0>0$ such that for any $0<\ga<\ga_0$, there holds $\lim_{\ge\to 0}A(\ge)=\infty$. This implies 
\begin{equation}\label{W5}
u_\infty(x(\ga),\ga)=\infty\qquad\forall \;0<\ga<\ga_0.
\end{equation}

\noindent{\it Step 2: Propagation.} In order to prove that the blow-up propagates along $\Gamma$ we have replace $t=0$ by $t=\ga<\ga_0$. We claim that
\begin{equation}\label{W6}\BA {ll}
\myint{B_\gs(x(\ga))}{}u_\infty(x,\ga)dx=\infty\qquad\forall \gs>0.
\EA\end{equation}
Actually, it is sufficient to prove the result with $\gs=\ge$ and with $u_\infty$ replaced by $v_\ge$. Then
\begin{equation}\label{W7}\BA {ll}
\myint{B_\ge(x(\ga)}{}v_\ge(x,\ga)dx=\myint{B_\ge}{}\tilde v_\ge(x,\ga)dx=\ge^{N-\frac{2}{p-1}}e^{\frac{\ell (\ge)}{p-1}}
\myint{B_1}{}\gw_\ge(x,\ge^{-2}\ga)dx
\EA\end{equation}
Using $(\ref{W0})$, $(\ref{W3})$ we have
$$\ge^{N-\frac{2}{p-1}}e^{\frac{\ell (\ge)}{p-1}}
\myint{B_1}{}\gw_\ge(x,\ge^{-2}\ga)dx\geq c'_1e^{A'(\ge)}\myint{B_1}{}\psi_0(x)dx
$$
where
\begin{equation}\label{W8}A'(\ge):=(N-\frac{2}{p-1})\ln\ge+\frac{\ell(\ge)}{p-1}-\left(\gl_0+\frac{\gb_\gt^2}{4}+\frac{\gd_\gt}{2}+\gs_\gt\right)\frac{\ga}{\ge^{2}}
\end{equation}
Thus $\lim_{\ge\to 0}A'(\ge)=\infty$ for any $\ga<\ga_0$. This implies the claim.
\smallskip

\noindent{\it Step 3: End of the proof.} 
For $k>0$, we denote by $u_{\ga,k\gd}$ the solution of 
\begin{equation}\label{Z1}\left\{\BA {ll}
\prt_tu-\Gd u+h(x,t)u^p=0\qquad&\text{in }\BBR^N\ti (\ga,\infty)\\
\phantom{\prt_tu+h(x,t)u^p}
u(.,\ga)=k\gd_{x(\ga)}\qquad&\text{in }\BBR^N.
\EA\right.\end{equation}
We claim that 
\begin{equation}\label{Z2}
u_{\infty}(x,t))\geq u_{\ga,k\gd}(x,t)\qquad\forall (x,t)\in \BBR^N\ti [\ga,\infty).
\end{equation}
We fix $k>0$, then for any $\gs>0$, there exists $m=m(\gs)>0$ such that 
$$\myint{B_\gs(x(\ga))}{}\min\{m,u_\infty(x,\ga)\}dx=k.
$$
Furthermore $\lim_{\gs\to 0}m(\gs)=\infty$. Let $u_{\ga,k,\gs}$ be the solution of 
\begin{equation}\label{Z3}\left\{\BA {ll}
\prt_tu-\Gd u+h(x,t)u^p=0\qquad&\text{in }\BBR^N\ti (\ga,\infty)\\
\phantom{\prt_tu+h(x,t)u^p}
u(.,\ga)=\min\{m,u_\infty(.,\ga)\}\chi_{_{B_\gs(x(\ga))}}\qquad&\text{in }\BBR^N.
\EA\right.\end{equation}
By the maximum principle $u_\infty\geq u_{\ga,k,\gs}$ in $\BBR^N\ti (\ga,\infty)$. But $\min\{m,u_\infty(.,\ga)\}\chi_{_{B_\gs(x(\ga))}}$ converges to $k\gd_{x(\ga)}$ in the weak sense of measure when $\gs\to 0$. By stability, since we have assumed $p<1+\frac{2}{N}$, $u_{\ga,k,\gs}\to u_{\ga,k\gd}$ locally uniformly in $\BBR^N\ti (\ga,\infty)$ (see \cite{MV2} e.g.). Therefore 
$(\ref{Z2})$ follows. Since $k$ is arbitrary, it follows
\begin{equation}\label{Z4}
u_{\infty}(x,t))\geq u_{\ga,\infty\gd}(x,t)\qquad\forall (x,t)\in \BBR^N\ti [\ga,\infty).
\end{equation}
By Step 1, $u_{\ga,\infty\gd}(x(\ga+\gb),\ga+\gb)$ is infinite for  $0<\gb<\ga_0$. This implies that $(\ref{W5})$ holds for any $0<\ga<2\ga_0$. Iterating this process we conclude that $u_\infty$ blows-up on whole $\Gamma$ and $(\ref{W6})$ holds for any $\ga\in (0,T]$.\qeda
\section{Propagation of singularities in the initial plane}
\setcounter{equation}{0}
In this section, we consider 
\begin{equation}\label{cur0}
\prt_{t}u-\Gd u+h(x,t)u^p=0\quad \text{in }Q_{\infty}:=\BBR^N\ti (0,\infty)
\end{equation}
where $h\in C(\overline Q_{\infty})$. We set $x=(x_1,...,x_N)=(x_1,x')$ and we suppose that $h(x,t)>0$ except when 
$(x,t)$ belongs to some straight line $\Gg$ that we can assume to be the $x_1$ axis in the plane $t=0$. We set
$d_\infty((x,t),\Gg)=\max\{\sqrt t,|x'|\}$ and write $h$ under the form
\begin{equation}\label{van1}
h(x,t)=e^{-\ell (d_\infty(x,t);\Gg)}
\end{equation}
where $\ell:(0,\infty)\mapsto (0,\infty)$ is continuous and nonincreasing with limit $\infty$ at $0$. \medskip

\noindent{\it Proof of Theorem D.} Case 1: $1<p<1+\frac{2}{N}$. For $\ge>0$, we consider the "tunnel" with axis $\Gamma$ defined by
$$T^\ge:=\{(x,t):x_1\in\BBR,(x',t)\in B'_\ge\ti (0,\ge^2)\}
$$
where $B'_\ge$ is the ball in $\BBR^{N-1}$ with center $0$ and radius $\ge$. Since $\ell$ is decreasing, there holds
$$
\prt_{t}u-\Gd u+e^{-\ell(\ge)}u^p\geq 0\quad \text{in }T_\ge.
$$
Thus 
$$u_{\infty\gd_0}\geq v_{\infty\gd_0}\quad \text{ in }T_\ge,
$$
where $v_{\infty\gd_0}=\lim_{k\to\infty}v_{k\gd_0}$ and $v_{k\gd}$ is the solution of
\begin{equation}\label{van2}\left\{\BA {ll}
\prt_{t}v-\Gd v+e^{-\ell(\ge)}v^p= 0\quad &\text{ in }T_\ge\\
\phantom{\prt_{t}-\Gd v+e^{-\ell(\ge)}v^p}
v=k\gd_0\quad &\text{ in }\BBR\ti B'_\ge\\
\phantom{\prt_{t}-\Gd v+e^{-\ell(\ge)}v^p}
v=0\quad &\text{ in }\BBR\ti \prt B'_\ge.
\EA\right.
\end{equation}
We put
$$v_{\infty\gd}(x,t)=\ge^{-2/(p-1)}e^{\ell(\ge)/(q-1)}w(x/\ge,t/\ge^2)
$$
Then $w=w_\ge$ satisfies
\begin{equation}\label{van3}\left\{\BA {ll}
\prt_{t}w-\Gd w+w^p= 0\quad &\text{ in }\BBR\ti B'_1\ti (0,1)\\
\phantom{\prt_{t}w-\Gd +w^p}
w=\infty\gd_0\quad &\text{ in }\BBR\ti B'_1\\
\phantom{\prt_{t}w-\Gd +w^p}
w=0\quad &\text{ in }\BBR\ti \prt B'_1.
\EA\right.
\end{equation}
We denote by $W$ the solution  of 
\begin{equation}\label{van4}\left\{\BA {ll}
\prt_{\gt}W-\Gd W+W= 0\quad &\text{ in }\BBR\ti B'_1\ti (0,1)\\
\phantom{\prt_{t,},-,}
W(\xi_1,\xi',0)=\psi(\xi_1)\phi(\xi')\quad &\text{ in }\BBR\ti B'_1\\
\phantom{\prt_{t,}W,-\Gd +W}
W=0\quad &\text{ in }\BBR\ti \prt B'_1\ti (0,1).
\EA\right.
\end{equation}
where $\phi$ is the first eigenfunction of $-\Gd_{\xi'}$ in $W^{1,2}_0(B_1')$ with maximum $1$ and corresponding eigenvalue $\gl$ and $\psi (\xi_1)=\cos (\xi_1)\chi_{_{[-\frac{\gp}{2},\frac{\gp}{2}]}}(\xi_1)$. Then
$$W(\xi_1,\xi',\gt)=\myfrac{e^{-(\gl+1) \gt}\phi(\xi')}{\sqrt{4\gp \gt}}\myint{-\gp/2}{\gp/2}e^{-|\xi_1-\gz|^2/4\gt}\psi (\gz)d\gz.
$$
Since $0\leq W\leq 1$,  $W$ is a subsolution for the equation satisfied by $w$ and there exists $a>0$ and $c>0$ such that
$$w_\ge(\xi,\gt+a)\geq cW(\xi,\gt)\qquad\forall (x,t)\in \BBR\ti B'_1\ti (0,1).
$$
Returning to $v_\infty\gd_0$, we derive 
\begin{equation}\label{van9}
v_{\infty\gd_0}(x,t+a\ge^2)\geq c\ge^{-2/(p-1)}e^{\ell(\ge)/(q-1)}W(x/\ge,t/\ge^2)\qquad\forall (x_1,x',t)\in \BBR\ti B'_\ge\ti(0,\ge^2],
\end{equation}
which implies, with $t=\ge^2$ and $x'=0$, 
\begin{equation}\label{van10}
v_\infty(x_1,0,(a+1)\ge^2)\geq c\ge^{-2/(p-1)}e^{\ell(\ge)/(q-1)}
\myfrac{e^{-\gl-1}\phi(0)}{\sqrt{4\gp}}\myint{-\gp/2}{\gp/2}e^{-|x_1/\ge-\gz|^2/4}\psi (\gz)d\gz.
\end{equation}
But
\begin{equation}\label{van11}
\myint{-\gp/2}{\gp/2}e^{-|x_1/\ge-\gz|^2/4}\psi (\gz)d\gz\geq 
e^{-x_1^2/2\ge^2}\myint{-\gp/2}{\gp/2}e^{-|\gz|^2/2}\psi (\gz)d\gz
\end{equation}
If we fix in particular $|x_1|\leq \gd$ where 
$$|x_1|<\gd=\sqrt\myfrac{2\ge^2\ell (\ge)}{q-1},
$$
we derive
\begin{equation}\label{van12}
\lim_{\ge\to 0}v_\infty(x_1,0,(a+1)\ge^2)=\infty.
\end{equation}
Furthermore this limit is uniform for $x_1\in [-\gd',\gd']$, where $\gd'<\gd$. Furthermore the interval  $[-\gd',\gd']$ does not shrink to $\{0\}$ when $\ge\to 0$, since it is assumed that
\begin{equation}\label{van13}
\liminf_{\ge\to 0}\ge^2\ell(\ge)>0.
\end{equation}
Replacing $[-\gd',\gd']$ by $[\gd',3\gd']$ and $[-3\gd',-\gd']$ and iterating, we conclude that
\begin{equation}\label{van14}
\lim_{t\to 0}u_{\infty\gd_0}(x_1,0,t)\geq\lim_{t\to 0}v_\infty(x_1,0,t)=\infty.
\end{equation}
From this, it is easy to obtain that $u_{\infty\gd_0}(x,t)=u_{\infty\gd_0}(0,x',t)$ is independent of $x_1$ and coincide with $U(x',t)$ where $U$ is the solution of 
\begin{equation}\label{van15}
\left\{\BA {ll}
\prt_{t}U-\Gd U+e^{-\ell(\max\{\sqrt t,|x'|\})}U^p= 0\quad &\text{ in } (0,\infty)\ti\BBR^{N-1}\\
\phantom{\prt_{t}U-\Gd +e^{-\ell(\max\{\sqrt t,|x'|\})}U^p}
U=\infty\gd_0\quad &\text{ in }\BBR^{N-1}\\
\phantom{\prt_{t}U-\Gd +U^p}
\EA\right.
\end{equation}

\noindent Case 2: $p\geq 1+\frac{2}{N}$. We write $h$ under the form
\begin{equation}\label{x1}
h(x,t)=\left(\max\{\sqrt t,|x'|\}\right)^{\gamma}e^{-\tilde \ell(\max\{\sqrt t,|x'|\})}
\end{equation}
where $\tilde \ell(s)=\ell(s)-\gamma\ln s$ and $\gamma>N(p-1)-2$. Then $\gamma>0$, $(\ref{I-3})$ is satisfied and for any $k>0$ there exists a unique solution $v=v_{k\gd}$ to
\begin{equation}\label{x2}\left\{\BA {ll}
\prt_{t}v-\Gd v+e^{-\tilde \ell(\ge)}\left(\max\{\sqrt t,|x'|\}\right)^{\gamma}v^p= 0\quad &\text{ in }T_\ge\\
\phantom{\prt_{t}v-\Gd +e^{-\tilde \ell(\ge)}\left(\max\{\sqrt t,|x'|\}\right)^{\gamma}v^p}
v=k\gd_0\quad &\text{ in }\BBR\ti B'_\ge\\
\phantom{\prt_{t}v-\Gd +e^{-\tilde \ell(\ge)}\left(\max\{\sqrt t,|x'|\}\right)^{\gamma}v^p}
v=0\quad &\text{ in }\BBR\ti \prt B'_\ge.
\EA\right.
\end{equation}
Furthermore $u_\infty(x,t)\geq v_{\infty\gd}$ in $T_\ge$. We set
$$v_{\infty\gd}(x,t)=\ge^{-(2+\gamma)/(p-1)}e^{\ell(\ge)/(q-1)}w(x/\ge,t/\ge^2),
$$
and $w=w_\ge$ satisfies

\begin{equation}\label{x3}\left\{\BA {ll}
\prt_{\gt}w-\Gd w+\left(\max\{\sqrt \gt,|\xi'|\}\right)^{\gamma}w^p= 0\quad &\text{ in }\BBR\ti B'_1\ti (0,1)\\
\phantom{\prt_{\gt}w-\Gd +\left(\max\{\sqrt \gt,|\xi'|\}\right)^{\gamma}w^p}
w=\infty\gd_0\quad &\text{ in }\BBR\ti B'_1\\
\phantom{\prt_{t}w-\Gd +\left(\max\{\sqrt \gt,|\xi'|\}\right)^{\gamma}w^p}
w=0\quad &\text{ in }\BBR\ti \prt B'_1\ti (0,1).
\EA\right.
\end{equation}
We denote by $W$ the solution of 
\begin{equation}\label{x4}\left\{\BA {ll}
\prt_{\gt}W-\Gd W+\left(\max\{\sqrt \gt,|\xi'|\}\right)^{\gamma}W= 0\quad &\text{ in }\BBR\ti B'_1\ti (0,1)\\
\phantom{,\prt_{t,},-\left(\max\{\sqrt \gt,|\xi'|\}\right)^{\gamma}}
W(\xi_1,\xi',0)=\psi(\xi_1)\phi(\xi')\quad &\text{ in }\BBR\ti B'_1\\
\phantom{\prt_{t,}W,-\Gd +W\left(\max\{\sqrt \gt,|\xi'|\}\right)^{\gamma}}
W=0\quad &\text{ in }\BBR\ti \prt B'_1\ti (0,1).
\EA\right.
\end{equation}
where $\psi$ and $\phi$ are as in the first case. This equation admits a separable solution $W(\xi,\gt)=W_1(\xi_1,\gt)W'(\xi',\gt)$ where
\begin{equation}\label{x5}\left\{\BA {ll}
\prt_{\gt}W'-\Gd_{\xi'} W'+\left(\max\{\sqrt \gt,|\xi'|\}\right)^{\gamma}W'= 0\quad &\text{ in }B'_1\ti (0,1)\\
\phantom{,-\Gd_{\xi'} W'+\left(\max\{\sqrt \gt,|\xi'|\}\right)^{\gamma}}
W'(\xi',0)=\phi(\xi')\quad &\text{ in } B'_1\\
\phantom{\prt_{\gt}W'-\Gd_{\xi'} W'+\left(\max\{\sqrt \gt,|\xi'|\}\right)^{\gamma}}
W'=0\quad &\text{ in }\prt B'_1\ti (0,1),
\EA\right.
\end{equation}
and
\begin{equation}\label{x6}\left\{\BA {ll}
\prt_{\gt}W_1-\prt_{x_1x_1}W_1+W_1= 0\quad &\text{ in }\BBR\ti (0,1)\\
\phantom{\prt_{\gt}W_1-\prt_{x_1x_1}}
W_1(\xi_1,0)=\phi(\xi_1)\quad &\text{ in } \BBR
\EA\right.
\end{equation}
Thus
$$
W(\xi_1,\xi',\gt)=\myfrac{W'(\xi',\gt)}{\sqrt{4\gp \gt}}\myint{-\gp/2}{\gp/2}e^{-|\xi_1-\gz|^2/4\gt}\psi (\gz)d\gz.
$$
The exact expression of $W'$ is not simple but since $\left(\max\{\sqrt \gt,|\xi'|\}\right)^{\gamma}\leq 1$ in $B'_1\ti (0,1)$, there holds
$$W(\xi_1,\xi',\gt)\geq \myfrac{e^{-(\gl+1) \gt}\phi(\xi')}{\sqrt{4\gp \gt}}\myint{-\gp/2}{\gp/2}e^{-|\xi_1-\gz|^2/4\gt}\psi (\gz)d\gz.
$$
From this point, and since $(\ref{van13})$ holds with $\ell$ replaced by $\tilde \ell$, the proof is the same as in  Case 1. \qeda

\end{document}